\documentclass[a4,12pt]{article}
\usepackage[centertags]{amsmath}
\usepackage{amsmath}
\usepackage{amssymb}
\usepackage{amsthm}
\usepackage{color}
\usepackage{dsfont}
\usepackage{latexsym}
\usepackage[T1]{fontenc}
\usepackage[latin1]{inputenc}
\usepackage{fancyhdr}
\usepackage[colorlinks=true,linkcolor=blue,citecolor=blue,urlcolor=black]{hyperref}
\newtheorem{definition}{Definition}
\newtheorem{theorem}{Theorem}
\newtheorem{proposition}{Proposition}
\newtheorem{lemma}{Lemma}

\newtheorem{corollary}{Corollary}

\def\ind{1\!\!1}
	\title{A Zero-Sum Deterministic Impulse Controls Game in Infinite Horizon with a New HJBI QVI}

	\author{Brahim El Asri\thanks{Ibn Zohr University, Equipe. Aide à la decision, ENSA, B.P. 1136, Agadir, Morocco. e-mail: b.elasri@uiz.ac.ma.} ,\, Hafid Lalioui\thanks{Ibn Zohr University, Equipe. Aide à la decision, ENSA, B.P. 1136, Agadir, Morocco. e-mail: hafid.lalioui@edu.uiz.ac.ma. Financially supported by CNRST, Rabat, Morocco (Grant 17 UIZ 2019).} \, and\, Sehail Mazid\thanks{Ibn Zohr University, Department of Mathematics, Faculty of Sciences, Morocco. e-mail: sehail.mazid@edu.uiz.ac.ma.}}

\begin{document}

\date{}
\maketitle

\begin{abstract}
	In the present paper, we study a two-player zero-sum deterministic differential game with both players adopting impulse controls, in infinite time horizon, under rather weak assumptions on the cost functions. We prove by means of the dynamic programming principle (DPP) that the lower and upper value functions are continuous and viscosity solutions to the corresponding Hamilton-Jacobi-Bellman-Isaacs (HJBI) quasi-variational inequality (QVI). We define a new HJBI QVI for which, under a proportional property assumption on the maximizer cost, the value functions are the unique viscosity solution. We then prove that the lower and upper value functions coincide.
\end{abstract}

\section{Introduction}
\quad

Differential games are concerned with the problem that multiple players make decisions, according to their own advantages and trade-off with other peers, in the context of dynamic systems. The theory of two-player zero-sum differential games was initiated by R. Isaacs \cite{I65} and L.S. Pontryagin et al \cite{PBGM62} at the beginning of 60's, in the early 80's the theory of viscosity solutions was pioneered by the seminal papers of M.G. Crandall and P.L. Lions \cite{CL83} and M.G. Crandall et al \cite{CEL84}. The notion of strategies and the rigorous definitions of lower and upper value functions are due to R.J. Elliott and N.J. Kalton \cite{EK72,EK74}, L.C. Evans and P.E. Souganidis \cite{ES84} began to study differential games by means of the viscosity theory, proving that the two value functions are the unique viscosity solution to the corresponding Hamilton-Jacobi-Bellman-Issacs (HJBI) partial differential equations (PDEs) for finite horizon problem, M. Bardi and I. Capuzzo-Dolcetta \cite{BC97} described the implementation of the viscosity solutions approach to a number of significant model problems in optimal deterministic control and differential games.

\par The deterministic differential games, apart from the mathematical interest in its own right, enjoy a wide range of applications in various fields of engineering, such as medicine, biology, economics and finance, see for more information A. Bensoussan and J.L. Lions \cite{BL84}. The deterministic impulse control problems in finite horizon were studied by many authors, J.M. Yong \cite{Y94} considered problems where one player takes continuous controls whereas the other uses impulse control, G. Barles et al \cite{ElBB10} treated a minimax problem driven by two controls, one is continuous and another impulsive. For the infinite horizon case as considered in the present paper, we cite the works of S. Dharmatti and A.J. Shaiju \cite{DS05,DS07}, S. Dharmatti and M. Ramaswamy \cite{DR06} and G. Barles \cite{B85}.

\par In this paper, we consider the state $y_{.}(.)$ of the two-player zero-sum deterministic differential game involving impulse controls in infinite time horizon described by the solution of the following system:
\begin{equation}\label{equation_1}
\left\{
\begin{aligned}
	y_x(0)&=x\in\mathbb{R}^n;\\
	\dot y_x(t)&=b\Bigl(y_x(t)\Bigr),\;t\neq\tau_m,\;t\neq\rho_k;\\
	y_x(\tau_m^+)&=y_x(\tau_m^-)+\xi_m\prod_{k\geq 1}\ind_{\{\tau_m\neq\rho_k\}},\;\tau_m\geq 0,\;\xi_m\neq 0;\\
	y_x(\rho_k^+)&=y_x(\rho_k^-)+\eta_k,\;\rho_k\geq 0,\;\eta_k\neq 0,
\end{aligned}
\right.
\end{equation}
where $m\in\mathbb{N^*}$, $k\in\mathbb{N^*}$ and $y_x(t)$ is the state variable of the system at time $t$, $\mathbb{R}^n$-valued, with initial state $y_x(0)=x$. The state $y_.(.)$ is driven by two impulse controls, $u$ control of $player-\xi$ defined by a double sequence $(\tau_m,\xi_m)_{m\geq 1}$ and $v$ control of $player-\eta$ defined by a double sequence $(\rho_k,\eta_k)_{k\geq 1}$, the actions $\xi_m$ and $\eta_k$ belong to the spaces of control actions $U\subset{\mathbb{R}^n}$ and $V\subset{\mathbb{R}^n}$, respectively. The infinite product $\prod_{k\geq 1}\ind_{\{\tau_m\neq\rho_k\}}$ signifies that when the two players act together on the system at the same time, we take into account only the action of $player-\eta$. The gain (resp. cost) functional $J$ for $player-\xi$ (resp. $player-\eta$) is defined as follows:
\begin{equation*}
J(x;u,v)=\int_0^\infty f\Bigl(y_x(t)\Bigr)exp(-\lambda t)dt-\sum_{m\geq 1}c(\xi_m)exp(-\lambda\tau_m)\prod_{k\geq 1}\ind_{\{\tau_m\neq\rho_k\}}+\sum_{k\geq 1}\chi(\eta_k)exp(-\lambda\rho_k),
\end{equation*}
where $c$ and $\chi$ are the zero lower bound impulse cost functions for $player-\xi$ and $player-\eta$, respectively, and the discount factor $\lambda$ is a positive real. We note that the cost of a player is the gain for the other (zero-sum), meaning that when a player performs an action he/she has to pay a positive cost, resulting in a gain for the other player. The function $f$ represents the running gain.

\par The terminology of a quasi-variational inequality (QVI) was introduced in \cite{BL84} to deal with impulse control problems. The definition of lower and upper value functions for a differential game as defined in \cite{EK72,EK74,ES84} leads to characterize the value of the game as the unique viscosity solution of a corresponding QVI. Moreover the relationship between the theory of two-player zero-sum deterministic differential games and viscosity solutions was first shown in \cite{ES84}, N. Barron et al \cite{BEJ84} and P.E. Souganidis \cite{S85,S85'}.

\par Recently, A. Cosso \cite{C13} and B. El Asri and S. Mazid \cite{ElM18} studied a stochastic impulse controls problem in finite horizon. They have proved, based on the dynamic programming principle (DPP) and viscosity solutions theory, that the differential game admits a value, however in both works the authors impose a stronger constraint that involves both cost functions, which is given by:
$$\exists h:[0,T]\rightarrow (0,+\infty)\;\text{such that}\; c(t,\xi_1+\eta+\xi_2)\leq c(t,\xi_1)-\chi(t,\eta)+c(t,\xi_2)-h(t),$$
where $c$ and $\chi$ are the impulse cost functions from $U$ and $V$, respectively, $\xi_1,\;\xi_2\in U$ and $\eta\in V$, as a consequence they had to require $V\subset U\subset\mathbb{R}^n$.

\par Our aim in this work is to investigate the two-player zero-sum deterministic impulse controls problem in infinite horizon given by the system (\ref{equation_1}). In particular, we describe the problem by a classic HJBI QVI, which we replace by a new HJBI QVI in order to characterize, in viscosity solution sense under rather weak assumptions on the cost functions, the value function of the differential game studied as the unique viscosity solution of the new HJBI QVI. In this work we only adopt, in addition to the classical assumptions of the impulse control problems, a proportional property assumption on the maximizer cost function $c$ which is given by:
\begin{equation}\label{equation_2}
\forall k>0,\;\forall\xi\in U\;\text{such that}\;k\xi\in U\;\text{we have}\;c(k\xi)\leq kc(\xi),
\end{equation}
note that assumption (\ref{equation_2}) is of great interest in the literature, as an application we cite the work developed in recent years in the field of biology, see L. Mailleret and F.Grognard \cite{MG09}.

\par For our game the associated QVI is given by the following double-obstacle HJBI equation, where the Hamiltonian involves only the first order partial derivatives:
\begin{equation}\label{equation_3}
max\Bigl\{min\Bigl[\lambda v(x)-Dv(x).b(x)-f(x),v(x)-\mathcal{H}_{sup}^c v(x)\Bigr],v(x)-\mathcal{H}_{inf}^\chi v(x)\Bigr\}=0,
\end{equation}
where $Dv(.)$ denotes the gradient of the function $v:\mathbb{R}^n\rightarrow\mathbb{R}$, and the first (resp. second) obstacle is defined through the use of the infimum (resp. maximum) cost operator $\mathcal{H}_{inf}^\chi$ (resp. $\mathcal{H}_{sup}^c$), where
$$\mathcal{H}_{inf}^\chi v(x)=\inf_{\eta\in V}\Bigl[v(x+\eta)+\chi(\eta)\Bigr]\;\Bigl(resp.\;\mathcal{H}_{sup}^c v(x)=\sup_{\xi\in U}\Bigl[v(x+\xi)-c(\xi)\Bigr]\Bigr).$$
We prove, using the dynamic programming principle, the existence of the value functions for our differential game as viscosity solutions of the HJBI QVI (\ref{equation_3}), but the uniqueness of the viscosity solution is not guarantee under standing assumptions, which means that the value function cannot enjoy anymore the property being the unique viscosity solution of the classic HJBI QVI (\ref{equation_3}).

\par Furthermore, we consider a new HJBI QVI, where the term of impulsions $v(.)-\mathcal{H}_{sup}^c v(.)$ (second obstacle) is replaced by the differential term $\mathcal F_{inf}^{c}\Bigl(Dv(.)\Bigr)$ defined for all $x\in\mathbb{R}^n$ by means of the operator $\mathcal F_{inf}^{c}$ as follows:
\begin{equation}\label{equation_4}
\mathcal F_{inf}^{c}\Bigl(Dv(x)\Bigr)=\inf_{\xi\in U}\Bigl[-Dv(x).\xi+c(\xi)\Bigr].
\end{equation}
Therefore, under assumption (\ref{equation_2}) and classical assumptions of the impulse control problems, we show the existence and the uniqueness results in the viscosity solution sense for the new HJBI QVI. Indeed, for the existence result, we give an equivalence in the viscosity supersolution sense between the classic HJBI QVI (\ref{equation_3}) and the new HJBI QVI, then, for the uniqueness, we establish a comparison theorem.

\par The paper is organized as follows: in section 2, we present the impulse controls problem studied, we give its related definitions and assumptions and we introduce our new HJBI QVI. In section 3, we prove some classic results on the lower and upper value functions, we first show that both satisfy the dynamic programming principle property, then we prove that they are continuous in $\mathbb{R}^n$. Section 4 is devoted, on one hand, to the viscosity characterization of the classic HJBI QVI (\ref{equation_3}) by deducing that both value functions are its viscosity solutions, and, on the other hand, we prove that the new HJBI QVI has the same bounded continuous viscosity supersolutions as the classic HJBI QVI (\ref{equation_3}), then we deduce the viscosity characterization of the new HJBI QVI. In the last section, we look more carefully to the new HJBI QVI by proving that the value functions of our infinite horizon two-player zero-sum impulse controls problem are his unique viscosity solution. Further, the lower and upper value functions coincide and the game admits a value.

\section{Assumptions and Setting of the Problem}
\subsection{Assumptions}
\quad

Throughout this paper, we let $n$ be a fixed positive integer, the time variable $T$ belongs to $[0,+\infty]$, $k$ and $m$ are in $\mathbb{N^*}$ and we let the discount factor $\lambda$ be a fixed positive real.
	
\par Let us assume \textbf{H1}:

\begin{itemize}

	\item[\textbf{[$H_{b,f}$]}]: The functions $b:\mathbb{R}^{n}\rightarrow\mathbb{R}^{n}$ and $f:\mathbb{R}^{n}\rightarrow\mathbb{R}$ are bounded and Lipschitz continuous with constant $C_b$ and $C_f$, respectively.
	\item[\textbf{[$H_{c,\chi}$]}]: The cost functions $c:U\rightarrow\mathbb{R}^{+}$ and $\chi:V\rightarrow\mathbb{R}^{+}$ are from two subsets of $\mathbb{R}^n$, $U$ and $V$, respectively, into $\mathbb{R}^{+}$, non negative and satisfy the zero lower bound property given by:
	$$\inf_{\xi\in U}c(\xi)> 0\;and\;\inf_{\eta\in V}\chi(\eta)>0.$$
	Also for all $\xi_1,\;\xi_2\in U$ and $\eta_1,\;\eta_2\in V$, we let the impulse cost functions satisfy
	$$c(\xi_1+\xi_2)\leq c(\xi_1)+c(\xi_2)$$
	and
	$$\chi(\eta_1+\eta_2)<\chi(\eta_1)+\chi(\eta_2).$$
	
\end{itemize}

\par Moreover, we assume \textbf{H2} that encompass the proportional impulse costs for $player-\xi$, that is the function $c$ satisfies
\begin{equation*}
\forall k>0,\;\forall\xi\in U\;\text{such that}\;k\xi\in U\;\text{we have}\;c(k\xi)\leq kc(\xi).
\end{equation*}

\par Regarding assumption \textbf{H1}, $[H_{b,f}]$ implies that there exists a unique global solution $y_{x}(.)$ to the above dynamical system (\ref{equation_1}), while the assumption $[H_{c,\chi}]$ provides the classical framework for the study of the impulse control problems. Assumption $\textbf{H2}$, which is of great interest in the literature, leads to the existence and the uniqueness results for the new HJBI QVI defined hereafter in (\ref{equation_6}).

\par For the rest of the paper we denote by $|.|$ and $\parallel.\parallel$ the Euclidian vector norm in $\mathbb{R}$ and $\mathbb{R}^n$, respectively, and for a bounded and continuous function $F$ from $\mathbb{R}^n$ to $\mathbb{R}$ (resp. $\mathbb{R}^n$) we define $\parallel F\parallel_\infty=\sup_{x\in\mathbb{R}^n}|F(x)|\;\Bigl(resp.\;\parallel F\parallel_\infty=\sup_{x\in\mathbb{R}^n}\parallel F(x)\parallel\Bigr).$

\subsection{Impulse Controls Game Problem}
\quad

Here we shall be interested in the two-player zero-sum deterministic differential game described in the introduction set by the dynamical system (\ref{equation_1}). The horizon (the interval in which time varies) is infinite. The state of the system $y_{.}(t)$ at the instant $t$ lies in $\mathbb{R}^n$, with initial value $y_x(0)=x$. The mapping $t\rightarrow y_x(t)$ describe the evolution of the system provided by a deterministic model $\dot y_x(t)=b\Bigl(y_x(t)\Bigr)$, where $b$ is a function from $\mathbb{R}^n$ to $\mathbb{R}^n$ satisfies assumption $[H_{b,f}]$.\\
At certain impulse instants $\tau_m$ and $\rho_k$, the state undergoes impulses (jumps) $\xi_m$ and $\eta_k$, respectively, that is:
\begin{equation*}
\begin{aligned}
	y_x(\tau_m^+)&=y_x(\tau_m^-)+\xi_m\prod_{k\geq 1}\ind_{\{\tau_m\neq\rho_k\}},\;\tau_m\geq 0,\;\xi_m\neq 0;\\
	y_x(\rho_k^+)&=y_x(\rho_k^-)+\eta_k,\;\rho_k\geq 0,\;\eta_k\neq 0,
\end{aligned}
\end{equation*}
the impulse time sequences $\{\tau_m\}_{m\geq 1}$ and $\{\rho_k\}_{k\geq 1}$ are two non-decreasing sequences of $[0,\infty]$ which satisfy $\tau_m,\;\rho_k\rightarrow+\infty$ when $m,\;k\rightarrow+\infty$, and the impulse value sequences $\{\xi_m\}_{m\geq 1}$ and $\{\eta_k\}_{k\geq 1}$ are two sequences of elements of $U\subset\mathbb{R}^n$ and $V\subset\mathbb{R}^n$, respectively.\\
The state $y_x(.)$ of the system is driven by two impulse controls, $(\tau_m,\xi_m)_{m\geq 1}$ control of $player-\xi$ and $(\rho_k,\eta_k)_{k\geq 1}$ control of $player-\eta$. The infinite product $\prod_{k\geq 1}\ind_{\{\tau_m\neq\rho_k\}}$ signifies that when the two players act together on the system at the same time, only the action of $player-\eta$ is tacking into account.

\par We call $\mathcal{U}$ (resp. $\mathcal{V}$) the space of impulse control $u$ (resp. $v$) for $player-\xi$ (resp. $player-\eta$) and we denote $u:=(\tau_m,\xi_m)_{m\geq 1}\;\Bigl(resp.\;v:=(\rho_k,\eta_k)_{k\geq 1}\Bigr)$. For any initial state $x$ the controls $u$ and $v$ generate a trajectory $y_x(.)$ solution of the system (\ref{equation_1}). We are given a gain (resp. cost) functional $J(x;u,v)$ for $player-\xi$ (resp. $player-\eta$), which represents the criterion to maximize (resp. minimize) by applying the control $u$ (resp. $v$):
\begin{equation*}
J(x;u,v)=\int_0^\infty f\Bigl(y_x(t)\Bigr)exp(-\lambda t)dt-\sum_{m\geq 1}c(\xi_m)exp(-\lambda\tau_m)\prod_{k\geq 1}\ind_{\{\tau_m\neq\rho_k\}}+\sum_{k\geq 1}\chi(\eta_k)exp(-\lambda\rho_k),
\end{equation*}
where $c$ and $\chi$ are the impulse cost functions (jump costs) from $U$ and $V$, respectively, $\mathbb{R}^+$-valued and satisfy assumption [$H_{c,\chi}$], the running gain $f:\mathbb{R}^n\rightarrow\mathbb{R}$ satisfies assumption [$H_{b,f}$] and $\lambda\geq 0$, representing the discount factor.

\par Typically, in a two-player game, the player who moves first would not choose a fixed action. Instead, he/she prefers to employ a strategy which can give different responses to different future actions the other player will take. Hence, besides the admissible controls, we define, following \cite{EK72,EK74}, the notion of non-anticipative strategies as follows:

\begin{definition}{(Non-anticipative strategy)}
	The non-anticipative strategy set $\mathcal{A}$ for $player-\xi$ is the collection of all non-anticipative maps $\alpha$ from $\mathcal{V}$ to $\mathcal{U}$, i.e., for any $v_1$ and $v_2$ in $\mathcal{V}$, if $v_1\equiv v_2$, then $\alpha(v_1)\equiv\alpha(v_2)$.\\
	The non-anticipative strategy set $\mathcal{B}$ for $player-\eta$ is the collection of all non-anticipative maps $\beta$ from $\mathcal{U}$ to $\mathcal{V}$, i.e., for any $u_1$ and $u_2$ in $\mathcal{U}$, if $u_1\equiv u_2$, then $\beta(u_1)\equiv\beta(u_2)$.
\end{definition}

In the game, $player-\xi$ aims to maximize the gain functional $J$ and contrarily $player-\eta$ aims to minimize. We may now give the definition of the lower and upper value functions for our two-player zero-sum deterministic differential game.\\ We define the \textit{lower value function} $V^-(.)$ and the \textit{upper value function} $V^+(.)$ of the game, respectively, by the following expressions:
\begin{equation*}
\begin{aligned}
	V^-(x)&=\inf_{\beta\in\mathcal B}\sup_{u\in\mathcal U} J\Bigl(x;u,\beta(u)\Bigr);\\
	V^+(x)&=\sup_{\alpha\in\mathcal A}\inf_{v\in\mathcal V} J\Bigl(x;\alpha(v),v\Bigr).
\end{aligned}
\end{equation*}
If $V^-=V^+$ we say that the game admits a value and $V:=V^-=V^+$ is called the \textit{value function} of the game.

\subsection{New HJBI Quasi-Variational Inequality}
\quad

For the impulse controls problem studied in the present paper, the associated Hamilton-Jacobi-Bellman-Isaacs quasi-variational inequality turns out to be the same for the two value functions, because of the two players can not act simultaneously on the system, and it is given by the following double-obstacle equation:
\begin{equation}\label{equation_5}
max\Bigl\{min\Bigl[\lambda v(x)-Dv(x).b(x)-f(x),v(x)-\mathcal{H}_{sup}^c v(x)\Bigr],v(x)-\mathcal{H}_{inf}^\chi v(x)\Bigr\}=0,
\end{equation}
where $\mathcal{H}_{inf}^\chi$ and $\mathcal{H}_{sup}^c$ are the cost operators defined, respectively, as follows:
\begin{equation*}
\begin{aligned}
	\mathcal{H}_{inf}^\chi v(x)&=\inf_{\eta\in V}\Bigl[v(x+\eta)+\chi(\eta)\Bigr];\\
	\mathcal{H}_{sup}^c v(x)&=\sup_{\xi\in U}\Bigl[v(x+\xi)-c(\xi)\Bigr],
\end{aligned}
\end{equation*}
and $Dv(.)$ denotes the gradient of the function $v:\mathbb{R}^n\rightarrow\mathbb{R}$.

\par Therefore, we define the new Hamilton-Jacobi-Bellman-Isaacs quasi-variational inequality (\ref{equation_6}), where the term of impulsions $v(.)-\mathcal{H}_{sup}^c v(.)$ is replaced by the differential term $\mathcal F_{inf}^{c}\Bigl(Dv(.)\Bigr)$, through the use of the operator $\mathcal F_{inf}^{c}$ as follows:
\begin{equation}\label{equation_6}
max\left\{min\left[\lambda v(x)-Dv(x).b(x)-f(x),\mathcal F_{inf}^{c}\Bigl(Dv(x)\Bigr)\right],v(x)-\mathcal{H}_{inf}^\chi v(x)\right\}=0,
\end{equation}
where the operator $\mathcal F_{inf}^{c}$ is defined as follows:
\begin{equation*}
\mathcal F_{inf}^{c}\Bigl(Dv(x)\Bigr)=\inf_{\xi\in U}\Bigl[-Dv(x).\xi+c(\xi)\Bigr].
\end{equation*}
Note that a differential term as in QVI (\ref{equation_6}) was introduced in G. Barles \cite{B94} and used, to deal with the particular case of null infimum jump costs in the infinite horizon impulse control problem, in N. El Farouq \cite{El17}.

\par The main objectives of this paper are:
\begin{itemize}
	\item[\textbf{\textit{(i)}}] Focusing on the existence of the solution in viscosity sense for both quasi-variational inequalities (\ref{equation_5}) and (\ref{equation_6}).
	\item[\textbf{\textit{(ii)}}] Showing that the QVI (\ref{equation_6}) admits the lower and upper value functions as the unique solution of viscosity.
\end{itemize}

\par For the rest of the paper we call QVI (\ref{equation_5}) the classic HJBI QVI, we call QVI (\ref{equation_6}) the new HJBI QVI and we adopt the following definition of the viscosity solution:

\begin{definition}{(Viscosity Solution)}
	Let $V:\mathbb{R}^n\rightarrow\mathbb{R}$ be a continuous function. $V$ is called:
	\begin{itemize}{}
		\item[\textbf{(i)}] A viscosity subsolution of the classic HJBI QVI (resp. new HJBI QVI) if for any $\overline x\in\mathbb{R}^n$ and any function $\phi\in C^{1}(\mathbb{R}^n)$ such that $V(\overline x)=\phi(\overline x)$ and $\overline x$ is a local maximum point of $V-\phi$, we have:
		$$max\Bigl\{min\Bigl[\lambda V(\overline x)-D\phi(\overline x).b(\overline x)-f(\overline x),V(\overline x)-\mathcal{H}_{sup}^c V(\overline x)\Bigr],V(\overline x)-\mathcal{H}_{inf}^\chi V(\overline x)\Bigr\}\leq 0,$$
		$$\Bigl(resp.\; max\Bigl\{min\Bigl[\lambda V(\overline x)-D\phi(\overline x).b(\overline x)-f(\overline x),\mathcal F_{inf}^c \Bigl(D\phi(\overline x)\Bigr)\Bigr],V(\overline x)-\mathcal{H}_{inf}^\chi V(\overline x)\Bigr\}\leq 0\Bigr).$$
		\item[\textbf{(ii)}] A viscosity supersolution of the classic HJBI QVI (resp. new HJBI QVI) if for any $\underline x\in\mathbb{R}^n$ and any function $\phi\in C^{1}(\mathbb{R}^n)$ such that $V(\underline x)=\phi(\underline x)$ and $\underline x$ is a local minimum point of $V-\phi$, we have:
		$$max\Bigl\{min\Bigl[\lambda V(\underline x)-D\phi(\underline x).b(\underline x)-f(\underline x),V(\underline x)-\mathcal{H}_{sup}^c V(\underline x)\Bigr],V(\underline x)-\mathcal{H}_{inf}^\chi V(\underline x)\Bigr\}\geq 0,$$
		$$\Bigl(resp.\; max\Bigl\{min\Bigl[\lambda V(\underline x)-D\phi(\underline x).b(\underline x)-f(\underline x),\mathcal F_{inf}^c \Bigl(D\phi(\underline x)\Bigr)\Bigr],V(\underline x)-\mathcal{H}_{inf}^\chi V(\underline x)\Bigr\}\geq 0\Bigr).$$
		\item[\textbf{(iii)}] A viscosity solution of the classic HJBI QVI (resp. new HJBI QVI) if it is both a viscosity subsolution and supersolution of the classic HJBI QVI (resp. new HJBI QVI).
	\end{itemize}
\end{definition}

\subsection{Preliminary Results}
\quad

Letting $y_x(.)$ and $y_{x^{'}}(.)$ be the trajectory generated by $u\in\mathcal{U}$ and $v:=\beta(u)\in\mathcal{V}$ from $x$ and $x^{'}$, respectively, where $\beta\in\mathcal{B}$. We then have the following characterization of the trajectories $y_{.}(.)$, for which the proof follows from Gronwall's Lemma and can be found in P.L. Lions \cite{L82}:
\begin{lemma}\label{lemma_1}
	Under assumption \textbf{H1}, for any $x,\;x^{'}\in\mathbb{R}^n$ and any $t\geq 0$ we have:
	$$\parallel y_x(t)-y_{x^{'}}(t)\parallel\leq exp(C_bt)\parallel x-x^{'}\parallel.$$
\end{lemma}
We show hereafter that the lower and upper value functions are bounded in $\mathbb{R}^n$.
\begin{proposition}\label{proposition_1}
	Under assumption \textbf{H1}, the lower and upper value functions are bounded in $\mathbb{R}^n$.
\end{proposition}
\begin{proof}
	We make the proof only for the lower value function $V^-$, the other case being analogous. By the definition of $V^-$, for all $x\in\mathbb{R}^n$ and all non-anticipative strategy $\beta\in\mathcal B$ we have
	$$V^-(x)\leq\sup_{u\in\mathcal U} J\Bigl(x;u,\beta(u)\Bigr).$$
	Then, considering the set of non-anticipative strategies $\beta(u):=(\rho_k,\eta_k)_{k\geq 1}$ for $player-\eta$ where there is no impulse time, i.e., $\rho_1=+\infty$, we get
	\begin{equation*}
	V^-(x)\leq\sup_{u\in\mathcal U}\biggl[\int_0^\infty f\Bigl(y_x(t)\Bigr)exp(-\lambda t)dt-\sum_{m\geq 1}c(\xi_m)exp(-\lambda\tau_m)\biggr].
	\end{equation*}
	Next, for all $\varepsilon>0$, there exists a strategy $u^\varepsilon:=(\tau_m^\varepsilon,\xi_m^\varepsilon)\in\mathcal{U}$ such that
	$$V^-(x)\leq\int_0^\infty f\Bigl(y_x(t)\Bigr)exp(-\lambda t)dt-\sum_{m\geq 1}c(\xi_m^\varepsilon)exp(-\lambda\tau_m^\varepsilon)+\varepsilon.$$
	Since $c$ is a non negative function and $f$ is bounded, we then get the existence of a constant $C>0$ such that
	$$V^-(x)\leq C.$$
	Similarly, for the set of strategies $u\in\mathcal{U}$ for $player-\xi$ for which there is no impulse time, i.e., $\tau_1=+\infty$, we have
	$$V^-(x)\geq\inf_{\beta\in\mathcal B}\biggl[\int_0^\infty f\Bigl(y_x(t)\Bigr)exp(-\lambda t)dt+\sum_{k\geq 1}\chi(\eta_k)exp(-\lambda\rho_k)\biggr].$$
	Let $\varepsilon>0$, then there exists a strategy $\beta^\varepsilon(u):=(\rho_k^\varepsilon,\eta_k^\varepsilon)\in\mathcal{V}$ where $\beta^\varepsilon\in\mathcal{B}$, for which we have
	$$V^-(x)+\varepsilon\geq\int_0^\infty f\Bigl(y_x(t)\Bigr)exp(-\lambda t)dt+\sum_{k\geq 1}\chi(\eta_k^\varepsilon)exp(-\lambda\rho_k^\varepsilon).$$
	Since $\chi$ is a non negative function and $f$ is bounded, we deduce the existence of a constant $C_1>0$ such that
	$$V^-(x)\geq -C_1.$$
	Hence we obtain the thesis.
\end{proof}

\section{Dynamic Programming Principle and Continuity of the Value Functions}
\quad

We now present, for our two-player zero-sum deterministic differential game in infinite time horizon, the dynamic programming principle in the following theorem. The DPP, as one of the principle and most commonly used approaches in solving optimal control problems, meaning that an optimal control viewed from today will remain optimal when viewed from tomorrow and stands for a basic property in dealing with our problem.

\begin{theorem}{(Dynamic Programming Principle)}
	Under assumption \textbf{H1}, given $x\in\mathbb{R}^n$ and $T>0$, we have the dynamic programming principle:
	\begin{equation}\label{equation_7}
	\begin{aligned}
		V^-(x)&=\inf_{\beta\in\mathcal B}\sup_{u\in\mathcal U}\biggl[\int_0^T f\Bigl(y_x(t)\Bigr)exp(-\lambda t)dt-\sum_{m\geq 1}c(\xi_m)exp(-\lambda\tau_m)\ind_{\{\tau_m\leq T\}}\prod_{k\geq 1}\ind_{\{\tau_m\neq\rho_k\}}\\
		&+\sum_{k\geq 1}\chi(\eta_k)exp(-\lambda\rho_k)\ind_{\{\rho_k\leq T\}}+V^-\Bigl(y_x(T)\Bigr)exp(-\lambda T)\biggr],
	\end{aligned}
	\end{equation}
	and
	\begin{equation*}
	\begin{aligned}
		V^+(x)&=\sup_{\alpha\in\mathcal A}\inf_{v\in\mathcal V}\biggl[\int_0^T f\Bigl(y_x(t)\Bigr)exp(-\lambda t)dt-\sum_{m\geq 1}c(\xi_m)exp(-\lambda\tau_m)\ind_{\{\tau_m\leq T\}}\prod_{k\geq 1}\ind_{\{\tau_m\neq\rho_k\}}\\
		&+\sum_{k\geq 1}\chi(\eta_k)exp(-\lambda\rho_k)\ind_{\{\rho_k\leq T\}}+V^+\Bigl(y_x(T)\Bigr)exp(-\lambda T)\biggr].
	\end{aligned}
	\end{equation*}
\end{theorem}

\begin{proof}
	We give the proof only for the lower value function $V^-$, similarly for $V^+$. We first let $\varepsilon>0$, $u\in\mathcal{U}$ and assume, for some $x\in\mathbb{R}^n$ and some $T>0$, that $V^-(x)<W_T(x)$, where $W_T(x)$ is the right-hand side of (\ref{equation_7}). We let the difference be $W_T(x)-V^-(x)=2\varepsilon$ and we choose $\beta^\varepsilon$ a non-anticipative strategy that approximates $V^-(x)$ up to $\varepsilon$, and denote $\beta^\varepsilon(u):=(\rho_k^\varepsilon,\eta_k^\varepsilon)_{k\geq 1}$ the jumps it produces. We then have
	\begin{equation*}
	\begin{aligned}
		W_T(x)-\varepsilon&\geq\sup_{u\in\mathcal U}\biggl[\int_0^T f\Bigl(y_x(t)\Bigr)exp(-\lambda t)dt-\sum_{m\geq 1}c(\xi_m)exp(-\lambda\tau_m)\ind_{\{\tau_m\leq T\}}\prod_{k\geq 1}\ind_{\{\tau_m\neq\rho_k^\varepsilon\}}\\
		&+\sum_{k\geq 1}\chi(\eta_k^\varepsilon)exp(-\lambda\rho_k^\varepsilon)\ind_{\{\rho_k^\varepsilon\leq T\}}+J\Bigl(y_x(T);u,\beta^\varepsilon(u)\Bigr)\biggr].
	\end{aligned}
	\end{equation*}
	The above inequality can be rewritten, for $\mathcal{U}_0$ and $\mathcal{U}_T$ the restrictions of $\mathcal{U}$ to $[0,T]\times U$ and $[T,+\infty]\times U$, respectively, as follows:
	\begin{equation}\label{equation_8}
	\begin{aligned}
		W_T(x)-\varepsilon&\geq\sup_{u\in\mathcal {U}_0}\biggl[\int_0^T f\Bigl(y_x(t)\Bigr)exp(-\lambda t)dt-\sum_{m\geq 1}c(\xi_m)exp(-\lambda\tau_m)\ind_{\{\tau_m\leq T\}}\prod_{k\geq 1}\ind_{\{\tau_m\neq\rho_k^\varepsilon\}}\\
		&+\sum_{k\geq 1}\chi(\eta_k^\varepsilon)exp(-\lambda\rho_k^\varepsilon)\ind_{\{\rho_k^\varepsilon\leq T\}}+\sup_{u\in\mathcal U_T}J\Bigl(y_x(T);u,\beta^\varepsilon(u)\Bigr)\biggr].
	\end{aligned}
	\end{equation}
	Observe that once $y_x(T)$ is known, the knowledge of $u$ over $[0,T]\times U$ is useless to evaluate $\sup_{u\in\mathcal U_T}J$. Therefore, since $\lambda\geq 0$, the restriction of $\mathcal{U}$ to $\mathcal U_T$ satisfies
	$$V^-\Bigl(y_x(T)\Bigr)exp(-\lambda T)\leq\sup_{u\in\mathcal U_T}J\Bigl(y_x(T);u,\beta^\varepsilon(u)\Bigr),$$
	replacing in inequality (\ref{equation_8}) leads to a contradiction. Finally, for all $x\in\mathbb{R}^n$ and all $T>0$, we deduce
	$$V^-(x)\geq W_T(x).$$
	Now let us assume to the contrary, for some $x\in\mathbb{R}^n$ and some $T>0$, that $V^-(x)>W_T(x)$ and let the difference be $V^-(x)-W_T(x)=3\varepsilon$ for $\varepsilon>0$. We denote by $\beta_1^\varepsilon\in\mathcal{B}$ the non-anticipative strategy that approximates $W_T(x)$ up to $\varepsilon$, where $\beta_1^\varepsilon(u):=(\rho_k^\varepsilon,\eta_k^\varepsilon)_{k\geq 1}$ are the jumps it produces for $u:=(\tau_m,\xi_m)_{m\geq 1}\in\mathcal{U}$. We then get
	\begin{equation*}
	\begin{aligned}
		V^-(x)-2\varepsilon&\geq\sup_{u\in\mathcal U}\biggl[\int_0^T f\Bigl(y_x(t)\Bigr)exp(-\lambda t)dt-\sum_{m\geq 1}c(\xi_m)exp(-\lambda\tau_m)\ind_{\{\tau_m\leq T\}}\prod_{k\geq 1}\ind_{\{\tau_m\neq\rho_k^\varepsilon\}}\\
		&+\sum_{k\geq 1}\chi(\eta_k^\varepsilon)exp(-\lambda\rho_k^\varepsilon)\ind_{\{\rho_k^\varepsilon\leq T\}}+V^-\Bigl(y_x(T)\Bigr)exp(-\lambda T)\biggr],
	\end{aligned}
	\end{equation*}
	therefore, for any $u\in\mathcal{U}_0$ the restriction of $\mathcal{U}$ to $[0,T]$, we have
	\begin{equation*}
	\begin{aligned}
		V^-(x)-2\varepsilon&\geq\int_0^T f\Bigl(y_x(t)\Bigr)exp(-\lambda t)dt-\sum_{m\geq 1}c(\xi_m)exp(-\lambda\tau_m)\ind_{\{\tau_m\leq T\}}\prod_{k\geq 1}\ind_{\{\tau_m\neq\rho_k^\varepsilon\}}\\
		&+\sum_{k\geq 1}\chi(\eta_k^\varepsilon)exp(-\lambda\rho_k^\varepsilon)\ind_{\{\rho_k^\varepsilon\leq T\}}+V^-\Bigl(y_x(T)\Bigr)exp(-\lambda T).
	\end{aligned}
	\end{equation*}
	Furthermore, we choose a non-anticipative strategy $\beta_2^\varepsilon\in\mathcal{B}$ of the game over $[T,+\infty]$ that approximates $V^-\Bigl(y_x(T)\Bigr)$ again up to $\varepsilon$. The concatenation $\beta^\varepsilon$ of $\beta_1^\varepsilon$ and $\beta_2^\varepsilon$ is a non-anticipative strategy of the game over $[0,+\infty]$, then we deduce for the non-anticipative strategy $\beta^\varepsilon\in\mathcal{B}$ and all control $u\in\mathcal{U}$ the following:
	$$V^-(x)-\varepsilon\geq J\Bigl(x;u,\beta^\varepsilon(u)\Bigr),$$
	a contradiction. Finally, for all $x\in\mathbb{R}^n$ and all $T>0$, we deduce
	$$V^-(x)\leq W_T(x).$$
	The proof is now complete.
\end{proof}

We next use the DPP to show that both the lower value function and the upper value function are continuous in $\mathbb{R}^n$.

\begin{theorem}
	Under assumption \textbf{H1}, the lower and upper value functions are continuous in $\mathbb{R}^n$.
\end{theorem}

\begin{proof}
	We make the proof only for the lower value function $V^-$, the other case being analogous. We first show that $V^-$ is upper semi-continuous. For any $x,\;x^{'}\in\mathbb{R}^n$ and $T>0$, according to the DPP for the lower value function, we have
	\begin{equation*}
	\begin{aligned}
		V^-(x)&=\inf_{\beta\in\mathcal{B}}\sup_{u\in\mathcal{U}}\biggl[\int_0^T f\Bigl(y_x(t)\Bigr)exp(-\lambda t)dt-\sum_{m\geq 1}c(\xi_m)exp(-\lambda\tau_m)\ind_{\{\tau_m\leq T\}}\prod_{k\geq 1}\ind_{\{\tau_m\neq\rho_k\}}\\
		&+\sum_{k\geq 1}\chi(\eta_k)exp(-\lambda\rho_k)\ind_{\{\rho_k\leq T\}}+V^-\Bigl(y_x(T)\Bigr)exp(-\lambda T)\biggr]
	\end{aligned}
	\end{equation*}
	and
	\begin{equation*}
	\begin{aligned}
		V^-(x^{'})&=\inf_{\beta\in\mathcal{B}}\sup_{u\in\mathcal{U}}\biggl[\int_0^T f\Bigl(y_{x^{'}}(t)\Bigr)exp(-\lambda t)dt-\sum_{m\geq 1}c(\xi_m)exp(-\lambda\tau_m)\ind_{\{\tau_m\leq T\}}\prod_{k\geq 1}\ind_{\{\tau_m\neq\rho_k\}}\\
		&+\sum_{k\geq 1}\chi(\eta_k)exp(-\lambda \rho_k)\ind_{\{\rho_k\leq T\}}+V^-\Bigl(y_{x^{'}}(T)\Bigr)exp(-\lambda T)\biggr].
	\end{aligned}
	\end{equation*}
	Now fix an arbitrary $\varepsilon>0$ and pick, for $\beta^\varepsilon\in\mathcal{B}$, a strategy $\beta^\varepsilon(u):=(\rho^\varepsilon_k,\eta^\varepsilon_k)_{k\geq 1}$ which satisfies the following:
	\begin{equation}\label{equation_9}
	\begin{aligned}
		V^-(x^{'})+\varepsilon&\geq\sup_{u\in\mathcal{U}}\biggl[\int_0^T f\Bigl(y_{x^{'}}(t)\Bigr)exp(-\lambda t)dt-\sum_{m\geq 1}c(\xi_m)exp(-\lambda\tau_m)\ind_{\{\tau_m\leq T\}}\prod_{k\geq 1}\ind_{\{\tau_m\neq\rho^\varepsilon_k\}}\\
		&+\sum_{k\geq 1}\chi(\eta^\varepsilon_k)exp(-\lambda\rho^\varepsilon_k)\ind_{\{\rho^\varepsilon_k\leq T\}}+V^-\Bigl(y_{x^{'}}(T)\Bigr)exp(-\lambda T)\biggr].
	\end{aligned}
	\end{equation}
	Further, we pick a control $u^\varepsilon:=(\tau^\varepsilon_m,\xi^\varepsilon_m)_{m\geq 1}\in\mathcal{U}$ which satisfies the following:
	\begin{equation*}
	\begin{aligned}
		V^-(x)-\varepsilon&\leq\int_0^T f\Bigl(y_x(t)\Bigr)exp(-\lambda t)dt-\sum_{m\geq 1}c(\xi^\varepsilon_m)exp(-\lambda\tau^\varepsilon_m)\ind_{\{\tau^\varepsilon_m\leq T\}}\prod_{k\geq 1}\ind_{\{\tau^\varepsilon_m\neq\rho^\varepsilon_k\}}\\
		&+\sum_{k\geq 1}\chi(\eta^\varepsilon_k)exp(-\lambda\rho^\varepsilon_k)\ind_{\{\rho^\varepsilon_k\leq T\}}+V^-\Bigl(y_x(T)\Bigr)exp(-\lambda T).
	\end{aligned}
	\end{equation*}
	Since the inequality (\ref{equation_9}) holds for any control $u$, then $u^\varepsilon$ satisfies
	\begin{equation*}
	\begin{aligned}
		V^-(x^{'})+\varepsilon&\geq\int_0^T f\Bigl(y_{x^{'}}(t)\Bigr)exp(-\lambda t)dt-\sum_{m\geq 1}c(\xi^\varepsilon_m)exp(-\lambda\tau^\varepsilon_m)\ind_{\{\tau^\varepsilon_m\leq T\}}\prod_{k\geq 1}\ind_{\{\tau^\varepsilon_m\neq\rho^\varepsilon_k\}}\\
		&+\sum_{k\geq 1}\chi(\eta^\varepsilon_k)exp(-\lambda\rho^\varepsilon_k)\ind_{\{\rho^\varepsilon_k\leq T\}}+V^-\Bigl(y_{x^{'}}(T)\Bigr)exp(-\lambda T).
	\end{aligned}
	\end{equation*}
	It follows from the two last inequalities that
	\begin{equation*}
	\begin{aligned}
		V^-(x)-V^-(x^{'})&\leq\int_0^T\Bigl[f\Bigl(y_x(t)\Bigr)-f\Bigl(y_{x^{'}}(t)\Bigr)\Bigr]exp(-\lambda t)dt\\
		&-\sum_{m\geq 1}c(\xi^\varepsilon_m)exp(-\lambda\tau^\varepsilon_m)\ind_{\{\tau^\varepsilon_m\leq T\}}\prod_{k\geq 1}\ind_{\{\tau^\varepsilon_m\neq\rho^\varepsilon_k\}}\\
		&+\sum_{k\geq 1}\chi(\eta^\varepsilon_k)exp(-\lambda\rho^\varepsilon_k)\ind_{\{\rho^\varepsilon_k\leq T\}}+V^-\Bigl(y_x(T)\Bigr)exp(-\lambda T)\\
		&+\sum_{m\geq 1}c(\xi^\varepsilon_m)exp(-\lambda\tau^\varepsilon_m)\ind_{\{\tau^\varepsilon_m\leq T\}}\prod_{k\geq 1}\ind_{\{\tau^\varepsilon_m\neq\rho^\varepsilon_k\}}\\
		&-\sum_{k\geq 1}\chi(\eta^\varepsilon_k)exp(-\lambda\rho^\varepsilon_k)\ind_{\{\rho^\varepsilon_k\leq T\}}-V^-\Bigl(y_{x^{'}}(T)\Bigr)exp(-\lambda T)+2\varepsilon.
	\end{aligned}
	\end{equation*}
	Thus, we get
	\begin{equation*}
	\begin{aligned}
		V^-(x)-V^-(x^{'})&\leq\int_0^T\parallel f\Bigl(y_x(t)\Bigr)-f\Bigl(y_{x^{'}}(t)\Bigr)\parallel exp(-\lambda t)dt\\
		&+\Bigl|V^-\Bigl(y_x(T)\Bigr)-V^-\Bigl(y_{x^{'}}(T)\Bigr)\Bigl|exp(-\lambda T)+2\varepsilon.
	\end{aligned}
	\end{equation*}
	By Lemma \ref{lemma_1}, the Lipschitz continuity of $f$ and the boundedness of $V^-$, we deduce that there exists a constant $C>0$ such that
	\begin{equation}\label{equation_10}
	V^-(x)-V^-(x^{'})\leq C_f\parallel x-x^{'}\parallel\int_0^T exp\Bigl((C_b-\lambda)t\Bigr)dt+2Cexp(-\lambda T)+2\varepsilon.
	\end{equation}
	Therefore, if $\lambda\neq C_b$, we obtain
	\begin{equation}\label{equation_11}
	V^-(x)-V^-(x^{'})\leq\frac{C_f}{C_b-\lambda}\parallel x-x^{'}\parallel\Bigl[exp\Bigl((C_b-\lambda)T\Bigr)-1\Bigr]+2Cexp(-\lambda T)+2\varepsilon.
	\end{equation}
	Now we choose $T$ such that $exp(-C_bT)=\parallel x-x^{'}\parallel^{1/2}$ with $\parallel x-x^{'}\parallel<1$. Hence, in the right-hand side of (\ref{equation_11}), the first term goes to $0$ when $x\rightarrow x^{'}\;, i.e.,\;T\rightarrow
	\infty$, indeed, it is equal to
	$$\frac{C_f}{C_b-\lambda}\parallel x-x^{'}\parallel^{1/2}\Bigl(exp(-\lambda T)-\parallel x-x^{'}\parallel^{1/2}\Bigr),$$
	while the second term goes to $0$ where $T\rightarrow\infty$. We then deduce, by letting $x\rightarrow x^{'}$ and $\varepsilon\rightarrow 0$, the upper semi-continuity of the lower value function:
	$$\limsup_{x\rightarrow x^{'}} V^-(x)\leq V^-(x^{'}).$$
	In the case where $\lambda=C_b$, it suffice to let some $\hat{\lambda}<\lambda=C_b$, so we go back to (\ref{equation_10}) and we proceed, since $exp\Bigl((C_b-\lambda) T\Bigr)<exp\Bigl((C_b-\hat{\lambda})T\Bigr)$ and $exp(-\lambda T)<exp(-\hat{\lambda}T)$, as above with the case $\hat{\lambda}\neq C_b$, we then conclude by letting $x\rightarrow x^{'}$ and $\varepsilon\rightarrow 0$.\\
	Analogously we get the lower semi-continuity:
	$$\liminf_{x\rightarrow x^{'}} V^-(x)\geq V^-(x^{'}).$$
	Then the lower value function is continuous in $\mathbb{R}^n$.
\end{proof}

\section{Viscosity Characterization}
\quad

The aim of the present section is to show that the lower and upper value functions are viscosity solutions to both classic HJBI QVI and new HJBI QVI. To do so we first give, in Lemma \ref{lemma_2}, some properties of the value functions, hence we get the aim for the classic HJBI QVI. Furthermore, we show the equivalence in viscosity supersolution sense between the classic HJBI QVI and the new HJBI QVI.\\
We begin with the following technical lemma:

\begin{lemma}\label{lemma_2}
	Under assumption \textbf{H1}, the lower value function satisfies for all $x\in\mathbb{R}^n$ the following property:
	$$V^-(x)\leq\mathcal{H}_{inf}^\chi V^-(x).$$
	Let $x\in\mathbb{R}^n$ be such that $V^-(x)<\mathcal{H}_{inf}^\chi V^-(x)$, then $V^-(x)\geq\mathcal{H}_{sup}^c V^-(x)$.\\
	The same results hold true for the upper value function $V^+$.
\end{lemma}

\begin{proof}
	We give the proof for the lower value function $V^-$, similarly for $V^+$. Letting $x\in\mathbb{R}^n$ and considering, for $player-\eta$, the strategy $\beta(u):=(\rho_k,\eta_k)_{k\geq 1}\in\mathcal{V}$ where $\beta\in\mathcal{B}$. Next, choose $\beta^{'}\in\mathcal{B}$ such that $\beta^{'}(u):=(0,\eta;\rho_2,\eta_2;\rho_3,\eta_3;.....)$, we then obtain $$V^-(x)\leq\sup_{u\in\mathcal{U}}J\Bigl(x;u,\beta^{'}(u)\Bigr)=\sup_{u\in\mathcal{U}}J\Bigl(x+\eta;u,\beta(u)\Bigr)+\chi(\eta),$$
	from which we deduce the following inequality:
	$$V^-(x)\leq\inf_{\eta\in\mathbb{R}^n}\Bigl[V^-(x+\eta)+\chi(\eta)\Bigr].$$
	Now let us assume that $V^-(x)<\mathcal{H}_{inf}^\chi V^-(x)$ for some $x\in\mathbb{R}^n$. From the DPP for $V^-$, by taking $T=0$, we get
	\begin{equation*}
	\begin{aligned}
		V^{-}(x)&=\inf_{\rho_1\in\{0,+\infty\},\;\eta\in V}\sup_{\tau_1\in\{0,+\infty\},\;\xi\in U}\biggl[ -c(\xi)\ind_{\{\tau_1=0\}}\ind_{\{\rho_1=+\infty\}}+\chi(\eta)\ind_{\{\rho_1=0\}}\\
		&+V^-(x+\xi\ind_{\{\tau_1=0\}}\ind_{\{\rho_1=+\infty\}}+\eta\ind_{\{\rho_1=0\}})\biggr],
	\end{aligned}
	\end{equation*}
	therefore
	\begin{equation*}
	\begin{aligned}
		V^-(x)&=\inf_{\rho_1\in\{0,+\infty\}}\Biggl[\inf_{\eta\in V}\Bigl[\chi(\eta)+V^-(x+\eta)\Bigr]\ind_{\{\rho_1=0\}}\\
		&+\sup_{\tau_1\in\{0,+\infty\},\;\xi\in U}\Bigl[-c(\xi)\ind_{\{\tau_1=0\}}+V^-(x+\xi\ind_{\{\tau_1=0\}})\Bigr]\ind_{\{\rho_1=+\infty\}}\Biggr].
	\end{aligned}
	\end{equation*}
	Since $V^-(x)<\mathcal{H}_{inf}^\chi V^-(x)$, we get
	$$V^-(x)=\sup_{\tau_1\in\{0,+\infty\},\;\xi\in U}\Bigl[-c(\xi)\ind_{\{\tau_1=0\}}+V^-(x+\xi\ind_{\{\tau_1=0\}})\Bigr].$$
	Therefore
	$$V^-(x)\geq\sup_{\xi\in U}\Bigl[V^-(x+\xi)-c(\xi)\Bigr].$$
\end{proof}

\subsection{Viscosity Characterization of the Classic HJBI Quasi-Variational Inequality}
\qquad

Now we are ready to show the relation between our deterministic impulse controls problem and the classic HJBI QVI (\ref{equation_5}), indeed, we prove the following theorem.

\begin{theorem}\label{theorem_3}
	Under assumption \textbf{H1}, the lower and upper value functions are viscosity solutions to the classic Hamilton-Jacobi-Bellman-Isaacs quasi-variational inequality.
\end{theorem}

\begin{proof}
	The proof is based on the DPP and it is inspired from \cite{B85}, we give the prove only for $V^-$, the other case being analogous. We first prove the subsolution property. Suppose $V^--\phi$ achieves its local maximum in $B_{\delta}(\overline x)$, where $B_{\delta}(\overline x)$ is the open ball of center $\overline x$ and radius $\delta>0$, with $V^-(\overline x)=\phi(\overline x)$, where $\phi$ is a function in $C^1(\mathbb{R}^n)$ and $x\in\mathbb{R}^n$. From Lemma \ref{lemma_2} we always have $V^-(\overline x)-\mathcal H_{inf}^\chi V^-(\overline x)\leq 0$, then if $V^-(\overline x)-\mathcal H_{sup}^c V^-(\overline x)\leq 0$ there is nothing to prove. Otherwise, we suppose that $V^-(\overline x)-\mathcal H_{sup}^c V^-(\overline x)>2\varepsilon_1>0$. Then, without loss of generality, we can assume that $V^-( x)-\mathcal H_{sup}^c V^-(x)>\varepsilon_1>0$ on $B_\delta(\overline x)$. Next, we define
	$$t'=\inf\Bigl\{t\geq 0: y_{\bar{x}}(t)\notin B_\delta(\overline x)\Bigr\}. $$
	Let $0<\varepsilon<\varepsilon_1$, $T\leq t'$ and consider $\rho_1=+\infty$, i.e., no impulse for $player-\eta$. Furthermore, we pick a strategy $u^\varepsilon:=(\tau^\varepsilon_m,\xi^\varepsilon_m)_{m\geq 1}\in\mathcal{U}$ for $player-\xi$ which, due to DPP for $V^-$, satisfies the following:
	\begin{equation*}
	\begin{aligned}
		V^-(\overline x)&\leq\int_0^{T\wedge\tau^\varepsilon_1} f\Bigl(y_{\overline x}(t)\Bigr)exp(-\lambda t)dt-c(\xi^\varepsilon_1)exp(-\lambda\tau^\varepsilon_1)\ind_{\{\tau^\varepsilon_1\leq T\}}\\
		&+V^-\Bigl(y_{\overline x}(T\wedge\tau^\varepsilon_1)\Bigr)exp\Bigl(-\lambda (T\wedge\tau^\varepsilon_1)\Bigr)+\varepsilon\\
		&\leq\int_0^{T\wedge\tau^\varepsilon_1} f\Bigl(y_{\overline x}(t)\Bigr)exp(-\lambda t)dt+exp(-\lambda \tau^\varepsilon_1)\biggl(V^-\Bigl(y_{\overline x}(\tau^\varepsilon_1)\Bigr)-c(\xi^\varepsilon_1)\biggr)\ind_{\{\tau^\varepsilon_1\leq T\}}\\
		&+V^-\Bigl(y_{\overline x}(T)\Bigr)exp(-\lambda T)\ind_{\{\tau^\varepsilon_1> T\}}+\varepsilon\\
		&\leq\int_0^{T\wedge\tau^\varepsilon_1} f\Bigl(y_{\overline x}(t)\Bigr)exp(-\lambda t)dt+\mathcal H_{sup}^c V^-\Bigl(y_{\overline x}(\tau^{\varepsilon-}_1)\Bigr)exp(-\lambda \tau^\varepsilon_1)\ind_{\{\tau^\varepsilon_1\leq T\}}\\
		&+V^-\Bigl(y_{\overline x}(T)\Bigr)exp(-\lambda T)\ind_{\{\tau^\varepsilon_1> T\}}+\varepsilon\\
		&\leq\int_0^{T\wedge\tau^\varepsilon_1} f\Bigl(y_{\overline x}(t)\Bigr)exp(-\lambda t)dt+V^-\Bigl(y_{\overline x}(\tau^{\varepsilon-}_1)\Bigr)exp(-\lambda \tau^\varepsilon_1)\ind_{\{\tau^\varepsilon_1\leq T\}}\\
		&+V^-\Bigl(y_{\overline x}(T)\Bigr)exp(-\lambda T)\ind_{\{\tau^\varepsilon_1> T\}}-\varepsilon_1 exp(-\lambda \tau^\varepsilon_1)\ind_{\{\tau^\varepsilon_1\leq T\}}+\varepsilon.\\
	\end{aligned}
	\end{equation*}
	Therefore, without loss of generality, we only need to consider strategy $u^\varepsilon\in\mathcal{U}$ such that $T<\tau^\varepsilon_1$, then we have
	\begin{equation*}
	V^-(\overline x)-\varepsilon\leq\int_0^{T} f\Bigl(y_{\overline x}(t)\Bigr)exp(-\lambda t)dt+V^-\Bigl(y_{\overline x}(T)\Bigr)exp(-\lambda T).
	\end{equation*}
	For $T$ small enough, we have
	$$\left\|y_{\overline x}(T)-\overline x\right\|\rightarrow 0,$$
	from which we deduce
	\begin{equation*}
	V^-\Bigl(y_{\overline x}(T)\Bigr)\leq\phi\Bigl(y_{\overline x}(T)\Bigr)+\Bigl[V^-(\overline x)-\phi(\overline x)\Bigr],
	\end{equation*}
	it follows, when $\varepsilon$ goes to $0$, that
	\begin{equation*}
	\dfrac{1-exp(-\lambda T)}{T}V^-(\overline x)\leq\dfrac{1}{T}\int_0^T f\Bigl(y_{\overline x}(t)\Bigr)exp(-\lambda t)dt+\dfrac{\phi\Bigl(y_{\overline x}(T)\Bigr)-\phi(\overline x)}{T}exp(-\lambda T).
	\end{equation*}
	We use the fact that
	$$\phi\Bigl(y_{\overline x}(T)\Bigr)-\phi(\overline x)=\int_{0}^{T}b\Bigl(y_{\overline x}(s)\Bigr).D\phi\Bigl(y_{\overline x}(s)\Bigr)ds,$$
	then we let $T\rightarrow 0$ to get
	$$\lambda V^-(\overline x)-D\phi(\overline x).b(\overline x)-f(\overline x)\leq 0.$$\\
	Finally, we get the subsolution property:
	$$max\Bigl\{min\Bigl[\lambda V^-(\overline x)-D\phi(\overline x).b(\overline x)-f(\overline x),V^-(\overline x)-\mathcal H_{sup}^c V^-(\overline x)\Bigr],V^-(\overline x)-\mathcal H_{inf}^\chi V^-(\overline x)\Bigr\}\leq 0.$$ The supersolution property is proved analogously.
\end{proof}

\subsection{Viscosity Characterization of the New HJBI Quasi-Variational Inequality}
\quad

This subsection is devoted to proving that both the lower and upper value functions of our problem are viscosity solutions to the new HJBI QVI, and we start by proving the following theorem.

\begin{theorem}\label{theorem_4}
	Under assumptions \textbf{H1} and \textbf{H2}, a bounded and continuous function $v$ is a viscosity supersolution to the classic HJBI QVI if and only if it is a viscosity supersolution to the new HJBI QVI.
\end{theorem}

\begin{proof}
	Assume first that $v$ is a bounded and continuous viscosity supersolution to the new HJBI QVI. Then, for a function $\phi\in C^{1}(\mathbb{R}^n)$ and $\underline x\in\mathbb{R}^n$ such that $v(\underline x)=\phi(\underline x)$ and $\underline x$ is a local minimum point of $v-\phi$, we have
	$$max\Bigl\{min\Bigl[\lambda v(\underline x)-D\phi(\underline x).b(\underline x)-f(\underline x),\mathcal F_{inf}^c D\phi(\underline x)\Bigr],v(\underline x)-\mathcal{H}_{inf}^\chi v(\underline x)\Bigr\}\geq 0.$$
	If $v(\underline x)-\mathcal{H}_{inf}^\chi v(\underline x)\geq 0$ then we are done with the first implication. Else we assume
	$$min\Bigl[\lambda v(\underline x)-D\phi(\underline x).b(\underline x)-f(\underline x),\mathcal F_{inf}^c D\phi(\underline x)\Bigr]\geq 0.$$
	It follows that $\mathcal F_{inf}^c D\phi(\underline x)\geq 0$. Let us first assume that $v\in C^{1}(\mathbb{R}^n)$, we then have for all $\xi\in U$, for all $x\in\mathbb{R}^n$,
	$$\mathcal F_{inf}^c Dv(x)\geq 0,$$
	that is, for all $\xi\in U$, for all $x\in\mathbb{R}^n$,
	$$-Dv(x).\xi\geq -c(\xi).$$
	Since
	$$v(x)-v(x+\xi)=\int_{0}^{1}-\frac{d}{ds}\Bigl(v(x+s\xi)\Bigr)ds=\int_{0}^{1}-Dv(x+s\xi).\xi ds\geq -c(\xi).$$
	We then obtain for all $\xi\in U$,
	$$v(\underline x)-\Bigl[v(\underline x+\xi)-c(\xi)\Bigr]\geq 0,$$
	which means that
	$$v(\underline x)-\mathcal{H}_{sup}^c v(\underline x)\geq 0.$$
	Finally, $v$ is a viscosity supersolution to the classic HJBI QVI when $v\in C^{1}(\mathbb{R}^n)$.\\
	We obtain the same result even $v$ is not in $C^{1}(\mathbb{R}^n)$. It suffices to make the same regularization $v_\varepsilon\in C^{\infty}(\mathbb{R}^n)$ for $v$, as in the proof of the Theorem $3.2$ in \cite{El17} (see also \cite{L82}), which converges uniformly toward $v$ in $\mathbb{R}^n$:\\
	Let $\theta$ be a positive function in $C^{\infty}(\mathbb{R}^n)$ with supp $\theta(x)\subset B_1(0)$, where $B_1(0)$ is the open ball of center $0$ and radius $1$, and $\int_{\mathbb{R}^n}\theta(x)dx=1$. We then define for $\varepsilon>0$ the function $\theta_\varepsilon$ by the following: $$\theta_\varepsilon(x)=\frac{1}{\varepsilon^n}\theta\Bigl(\frac{x}{\varepsilon}\Bigl).$$
	Further, we define in $\mathbb{R}^n$ the regularization
	$$v_\varepsilon(x)=\int_{\mathbb{R}^n}v(y)\theta_\varepsilon(x-y)dy.$$
	The function $v_\varepsilon$ is bounded, belongs to $C^{\infty}(\mathbb{R}^n)$ and satisfies
	$$\sup_{x\in\mathbb{R}^n}\Bigl|v_\varepsilon(x)-v(x)\Bigr|\leq\sup_{\parallel x-y\parallel\leq\varepsilon}\Bigl|v(x)-v(y)\Bigr|,$$
	then it converges to $v$ as $\varepsilon$ goes to $0$.\\
	In addition, for all $x\in\mathbb{R}^n$ and since $v$ is a viscosity supersolution to the new HJBI QVI, the regularization $v_\varepsilon$ satisfies $\mathcal F_{inf}^c Dv_\varepsilon(x)\geq\delta(\varepsilon)$, where $\delta(\varepsilon)$ goes to $0$ with $\varepsilon$.\\
	Therefore, by the same computation as in above, we get for all $\xi\in U$,
	$$v_\varepsilon(\underline x)-\Bigl[v_\varepsilon(\underline x+\xi)-c(\xi)\Bigr]\geq\delta(\varepsilon),$$
	that is, when $\varepsilon$ goes to $0$,
	$$v(\underline x)-\Bigl[v(\underline x+\xi)-c(\xi)\Bigr]\geq 0.$$
	Hence we get the desired result:
	$$v(\underline x)-\mathcal{H}_{sup}^c v(\underline x)\geq 0.$$
	Assume now that $v$ is a bounded and continuous viscosity supersolution to the classic HJBI QVI. And let $\underline x\in\mathbb{R}^n$ be a global minimum point of $v-\phi$, where $\phi$ is a function in $C^{1}(\mathbb{R}^n)$ and $v(\underline x)=\phi(\underline x)$. We then have
	$$max\Bigl\{min\Bigl[\lambda v(\underline x)-D\phi(\underline x).b(\underline x)-f(\underline x),v(\underline x)-\mathcal{H}_{sup}^c v(\underline x)\Bigr],v(\underline x)-\mathcal{H}_{inf}^\chi v(\underline x)\Bigr\}\geq 0.$$
	If $v(\underline x)-\mathcal{H}_{inf}^\chi v(\underline x)\geq 0$ then we are done. Else we assume
	$$min\Bigl[\lambda v(\underline x)-D\phi(\underline x).b(\underline x)-f(\underline x),v(\underline x)-\mathcal{H}_{sup}^c v(\underline x)\Bigr]\geq 0.$$
	Which gives
	$$v(\underline x)-\sup_{\xi\in U}\Bigl[v(\underline x+\xi)-c(\xi)\Bigr]\geq 0,$$
	thus, for all $\xi\in U$, we get
	$$v(\underline x)-v(\underline x+\xi)+c(\xi)\geq 0.$$
	Since for all $\xi\in U$,
	$$\phi(\underline x)-\phi(\underline x+\xi)+c(\xi)\geq v(\underline x)-v(\underline x+\xi)+c(\xi),$$
	then we have
	$$\phi(\underline x)-\phi(\underline x+\xi)+c(\xi)\geq 0.$$
	We can then deduce, under assumption \textbf{H2}, for all $k>0$ the following:
	$$\frac{\phi(\underline x)-\phi(\underline x+k\xi)}{k}\geq\frac{-c(k\xi)}{k}\geq -c(\xi).$$
	Hence, by letting $k\rightarrow 0$, we get for all $\xi\in U$,
	$$-D\phi(\underline x).\xi+c(\xi)\geq 0.$$
	Therefore, we obtain
	$$\mathcal F_{inf}^c D\phi(\underline x)\geq 0.$$
	Finally, $v$ is a viscosity supersolution to the new HJBI QVI.
\end{proof}

\begin{theorem}
	Under assumptions \textbf{H1} and \textbf{H2}, the lower and upper value functions are viscosity solutions to the new Hamilton-Jacobi-Bellman-Isaacs quasi-variational inequality.
\end{theorem}

\begin{proof}
	We give the proof for the lower value function $V^-$, similarly for $V^+$. We first prove the viscosity subsolution property. We let $\phi$ be a function in $C^{1}(\mathbb{R}^n)$ and $\overline x\in\mathbb{R}^n$ such that $V^-(\overline x)=\phi(\overline x)$ and $\overline x$ is a local maximum of $V^{-}-\phi$. Since we have proved in Lemma \ref{lemma_2} that $V^-(\overline x)-\mathcal{H}_{inf}^\chi V^-(\overline x)\leq 0$, then if $\mathcal F_{inf}^c D\phi(\overline x)\leq 0$ there is nothing to prove. Otherwise, for all $x\in\mathbb{R}^n$ we assume that $\mathcal F_{inf}^c D\phi(x)>0$. We then get for all $\xi\in U$,
	$$\inf_{\xi\in U}\Bigl[-D\phi(x).\xi+c(\xi)\Bigr]>0.$$
	Then for all $\xi\in U$, for all $x\in\mathbb{R}^n$,
	$$-D\phi(x).\xi>-c(\xi).$$
	Since
	$$\phi(x)-\phi(x+\xi)=\int_{0}^{1}-\frac{d}{ds}\Bigl(\phi(x+s\xi)\Bigr)ds=\int_{0}^{1}-D\phi(x+s\xi).\xi ds>-c(\xi),$$
	we get
	$$\phi(\overline x)-\phi(\overline x+\xi)+c(\xi)>0,$$
	furthermore
	$$\phi(\overline x)-\phi(\overline x+\xi)+c(\xi)\leq V^-(\overline x)-V^-(\overline x+\xi)+c(\xi).$$
	Hence
	$$V^-(\overline x)-\Bigl[V^-(\overline x+\xi)-c(\xi)\Bigr]> 0.$$
	Thus, whenever $\mathcal F_{inf}^c D\phi(\overline x)>0$ we have $V^-(\overline x)-\mathcal{H}_{sup}^c V^-(\overline x)>0$.\\
	Next, the same computation as in the proof of the viscosity subsolution sense for the classic HJBI QVI, Theorem \ref{theorem_3}, leads to the following viscosity subsolution property:
	$$max\Bigl\{min\Bigl[\lambda V^-(\overline x)-D\phi(\overline x).b(\overline x)-f(\overline x),\mathcal F_{inf}^c D\phi(\overline x)\Bigr],V^-(\overline x)-\mathcal{H}_{inf}^\chi V^-(\overline x)\Bigr\}\leq 0.$$
	For the prove of the viscosity supersolution property, since we have proved in Theorem \ref{theorem_3} that $V^-$ is a viscosity solution to the classic HJBI QVI, then, due to the result in Theorem \ref{theorem_4}, the lower value function $V^-$ is a viscosity supersolution to the new HJBI QVI.
\end{proof}

\section{Uniqueness of the Viscosity Solution of the New HJBI Quasi-Variational Inequality}
\qquad

We prove in the present section, via a comparison theorem, that the new HJBI QVI has a unique bounded and continuous solution in viscosity sense. As a consequence, the lower and upper value functions coincide, since they are both viscosity solutions to the new HJBI QVI. Hence the game has a value.\\
We first give the following classical lemma which also appears in \cite{Y94}:

\begin{lemma}\label{lemma_3}
	Let $v$ be a bounded uniformly continuous function and $x_0\in\mathbb{R}^n$ such that:
	$$v(x_0)\geq\mathcal{H}_{inf}^\chi v(x_0).$$
	Then there exists an element $y$ in $\mathbb{R}^n$ such that:
	$$\exists\delta>0,\;\forall x\in\overline{B}_\delta(y):\;v(x)<\mathcal{H}_{inf}^\chi v(x),$$
	where $\overline{B}_\delta(y)$ is the closed ball of center $y$ and radius $\delta$.
\end{lemma}

\begin{proof}
	Fix $\varepsilon>0$ and let $x_0\in\mathbb{R}^n$ such that $v(x_0)\geq\mathcal{H}_{inf}^\chi v(x_0)$. Then there exists $\eta_0\in V$ such that
	$$v(x_0)\geq v(x_0+\eta_0)+\chi(\eta_0)-\varepsilon,$$
	then we get for all $\eta\in V$,
	\begin{equation*}
	\begin{aligned}
		v(x_0+\eta_0+\eta)+\chi(\eta)-v(x_0+\eta_0)&\geq v(x_0+\eta_0+\eta)+\chi(\eta)+\chi(\eta_0)-v(x_0)-\varepsilon\\
		&\geq\chi(\eta_0)+\chi(\eta)-\chi(\eta_0+\eta)-\varepsilon.
	\end{aligned}
	\end{equation*}
	Thus, from assumption \textbf{H1} and by letting $\varepsilon\rightarrow 0$, we deduce
	$$v(x_0+\eta_0)<\mathcal{H}_{inf}^\chi v(x_0+\eta_0).$$
	Now we take $y=x_0+\eta_0$, then, since for all $\eta\in V$ we have
	$$v(y)<v(y+\eta)+\chi(\eta)-\varepsilon,$$
	we obtain, by uniform continuity of $v$ where $C_v$ is the modulus of continuity, for all $x\in\mathbb{R}^n$,
	$$v(x)<v(x+\eta)+\chi(\eta)+C_v(\parallel x-y\parallel)-\varepsilon.$$
	Hence there exists $\delta>0$ such that for all $x\in\overline{B}_\delta(y)$ we have
	$$v(x)-\mathcal{H}_{inf}^\chi v(x)<0.$$
\end{proof}

We next prove the following useful lemma.

\begin{lemma}\label{lemma_4}
	Let $v:\mathbb{R}^n\rightarrow\mathbb{R}$ be a bounded and continuous viscosity supersolution to the new HJBI QVI. If $v(x)-\mathcal{H}_{inf}^\chi v(x)<0$, then, for any $\mu$, $\alpha$ and $K$ such that
	$$0<\mu<1,\;K>\parallel f\parallel_\infty/\lambda\;\text{and}\;0<\alpha<(1-\mu)\min\Bigl(\inf_{\xi\in U}c(\xi),(\lambda K-\parallel f\parallel_\infty)/\parallel b\parallel_\infty\Bigr),$$
	the function $v^*(x)=\mu v(x)+\alpha\sqrt{\parallel x\parallel^2+1}+K(1-\mu)$ is a strict viscosity supersolution to the new HJBI QVI.
\end{lemma}

\begin{proof}
	Let $v$ be a bounded continuous viscosity supersolution to the new HJBI QVI, $\phi^*\in C^{1}(\mathbb{R}^n)$ and $\underline x\in\mathbb{R}^n$ be a local minimum point of $v^*-\phi^*$ such that $v(\underline x)-\mathcal{H}_{inf}^\chi v(\underline x)<0$.\\
	Then, for all $x\in B_\delta(\underline x)$, where $B_\delta(\underline x)$ is the open ball of center $\underline x$ and radius $\delta>0$, we have
	$$\mu v(\underline x)+\alpha\sqrt{\parallel\underline x\parallel^2+1}-\phi^*(\underline x)\leq\mu v(x)+\alpha\sqrt{\parallel x\parallel^2+1}-\phi^*(x),$$
	then
	$$v(\underline x)-\frac{\phi^*(\underline x)-\alpha\sqrt{\parallel\underline x\parallel^2+1}}{\mu}\leq v(x)-\frac{\phi^*(x)-\alpha\sqrt{\parallel x\parallel^2+1}}{\mu},$$
	this inequality means that $\underline x$ is a local minimum point of $v-\phi$, where $$\phi(x)=\Bigl(\phi^*(x)-\alpha\sqrt{\parallel x\parallel^2+1}\Bigr)/\mu.$$
	Then, since $v(\underline x)-\mathcal{H}_{inf}^\chi v(\underline x)<0$ and $v$ is viscosity supersolution to the new HJBI QVI, we have the following:
	$$min\Bigl[\lambda v(\underline x)-D\phi(\underline x).b(\underline x)-f(\underline x),\mathcal F_{inf}^c D\phi(\underline x)\Bigr]\geq 0,$$
	thus
	$$\lambda v(\underline x)-D\phi(\underline x).b(\underline x)-f(\underline x)\geq 0\;\text{and}\;\mathcal F_{inf}^c D\phi(\underline x)\geq 0.$$
	On one hand, we get
	$$\lambda v(\underline x)-\frac{1}{\mu}D\phi^*(\underline x).b(\underline x)+\frac{\alpha\underline x}{\mu\sqrt{\parallel\underline x\parallel^2+1}}.b(\underline x)-f(\underline x)\geq 0,$$
	then
	$$\lambda\mu v(\underline x)-D\phi^*(\underline x).b(\underline x)-\mu f(\underline x)\geq-\alpha\parallel b\parallel_\infty,$$
	thus
	$$\lambda v^*(\underline x)-D\phi^*(\underline x).b(\underline x)-f(\underline x)\geq \lambda K(1-\mu)-(1-\mu)\parallel f\parallel_\infty-\alpha\parallel b\parallel_\infty.$$
	Finally, since $\alpha<(1-\mu)(\lambda K-\parallel f\parallel_\infty)/\parallel b\parallel_\infty$, we obtain
	\begin{equation}\label{equation_12}
	\lambda v^*(\underline x)-D\phi^*(\underline x).b(\underline x)-f(\underline x)> 0.
	\end{equation}
	On the other hand, since $\mathcal F_{inf}^c D\phi(\underline x)\geq 0,$ we get
	$$\mathcal F_{inf}^c\biggl[\frac{D\phi^*(\underline x)}{\mu}-\frac{\alpha\underline x}{\mu\sqrt{\parallel\underline x\parallel^2+1}}\biggr]\geq 0,$$
	then
	$$\inf_{\xi\in U}\biggl[-D\phi^*(\underline x).\xi+\frac{\alpha\underline x}{\sqrt{\parallel\underline x\parallel^2+1}}.\xi+\mu c(\xi)\biggr]\geq 0,$$
	it follows that
	$$\inf_{\xi\in U}\Bigl[-D\phi^*(\underline x).\xi+c(\xi)+(\mu-1)c(\xi)\Bigr]\geq -\alpha,$$
	from which we deduce, since $\mu-1<0$, the following:
	$$\inf_{\xi\in U}\Bigl[-D\phi^*(\underline x).\xi+c(\xi)+(\mu-1)\inf_{\xi\in U}c(\xi)\Bigr]\geq\inf_{\xi\in U}\Bigl[-D\phi^*(\underline x).\xi+c(\xi)+(\mu-1)c(\xi)\Bigr]\geq-\alpha,$$
	therefore
	$$(\mu-1)\inf_{\xi\in U}c(\xi)+\inf_{\xi\in U}\Bigl[-D\phi^*(\underline x).\xi+c(\xi)\Bigr]\geq-\alpha.$$
	Finally, from assumption \textbf{H1} and since $\alpha<(1-\mu)\inf_{\xi\in U}c(\xi)$, we obtain
	\begin{equation}\label{equation_13}
	\mathcal F_{inf}^c D\phi^*(\underline x)>0.
	\end{equation}
	Then the two strict inequalities (\ref{equation_12}) and (\ref{equation_13}) imply that
	$$min\Bigl[\lambda v^*(\underline x)-D\phi^*(\underline x).b(\underline x)-f(\underline x),\mathcal F_{inf}^c D\phi^*(\underline x)\Bigr]>0.$$
	Thus, we get that $v^*$ is a strict viscosity supersolution to the new HJBI QVI.
\end{proof}

We are now in a position to prove the comparison theorem.

\begin{theorem}{(Comparison Theorem)}
	Under assumptions \textbf{H1} and \textbf{H2}, if $u$ is a bounded and continuous viscosity subsolution to the new HJBI QVI and $v$ is a bounded and continuous viscosity supersolution to the new HJBI QVI, then we have: $$\forall x\in\mathbb{R}^n:\;u(x)\leq v(x).$$
\end{theorem}

\begin{proof}
	Let $u$ and $v$ be a bounded and continuous viscosity subsolution and supersolution, respectively, to the new HJBI QVI.
	Our aim is to show, by contradiction, that $u\leq v$.\\
	We denote by $M=\sup_{x\in\mathbb{R}^n}\Bigl(u(x)-v^*(x)\Bigr)$ the maximal value of $u-v^*$, where $v^*$ is defined as in Lemma \ref{lemma_4}. We let $R=\Bigl(\parallel u\parallel_\infty+\parallel v\parallel_\infty\Bigr)/\alpha$, then we have for all $x\in\mathbb{R}^n$ such that $\parallel x\parallel\geq R$,
	$$u(x)\leq\parallel u\parallel_\infty+(1-\mu)\parallel v\parallel_\infty\leq v^*(x),$$
	that is $u(x)\leq v^*(x)$ for all $x\in\mathbb{R}^n\textbackslash \overline{B}_{R}(0)$ where $\overline{B}_{R}(0)$ is the closed ball in $\mathbb{R}^n$ of radius $R$ centered at $0$.\\
	Let us now assume that there exists $\hat{x}\in B_{R}(0)$, the open ball, such that
	$$M=u(\hat{x})-v^*(\hat{x})>0,$$
	if it is not the case, i.e., $M\leq 0$, then the proof is finished, will follow by letting $\mu\rightarrow 1$ and $\alpha\rightarrow 0$.
	\par\textbf{Step 1.}
	We can find $\overline{x}\in B_{R}(0)$ and $\delta>0$ such that
	$$\sup_{x\in\overline{B}_\delta(\overline{x})}\Bigl(u(x)-v^*(x)\Bigr)\geq u(\overline{x})-v^*(\overline{x})>0,$$
	and for all $x\in\overline{B}_\delta(\overline{x}),$
	$$v(x)<\mathcal{H}_{inf}^\chi v(x),$$
	where $\overline{B}_\delta(\overline{x})$ is the closed ball of center $\overline{x}$ and radius $\delta$.\\
	In fact, if $v(\hat{x})<\mathcal{H}_{inf}^\chi v(\hat{x})$, then considering the continuity of $u$, $v$ and $\mathcal{H}_{inf}^\chi$ we obtain the result by taking $\overline{x}=\hat{x}$.\\
	Otherwise, we let $v(\hat{x})\geq\mathcal{H}_{inf}^\chi v(\hat{x})$, then, for some $\eta^{'}\in V$, the result in Lemma \ref{lemma_3} gives
	$$v(\hat{x}+\eta^{'})<\mathcal{H}_{inf}^\chi v(\hat{x}+\eta^{'}),$$
	we then take $\overline{x}=\hat{x}+\eta^{'}$ to deduce
	\begin{equation}\label{equation_14}
	\exists\delta>0,\;\forall x\in\overline{B}_\delta(\overline{x}):\;v(x)<\mathcal{H}_{inf}^\chi v(x).
	\end{equation}
	Furthermore, when $v^*(\hat{x})\geq\mathcal{H}_{inf}^\chi v^*(\hat{x})$, we fix $\varepsilon>0$, then for $\alpha\in(0,1)$ there exists $\eta^{'}\in V$ such that
	\begin{equation}\label{equation_15}
	v^*(\hat x)\geq v^*(\hat x+\eta^{'})+\chi(\eta^{'})-\alpha\varepsilon,
	\end{equation}
	which gives
	\begin{equation*}
	\begin{aligned}
		u(\hat x+\eta^{'})-v^*(\hat x+\eta^{'})&\geq u(\hat x+\eta^{'})+\chi(\eta^{'})-v^*(\hat x)-\alpha\varepsilon\\
		&\geq u(\hat x)-v^*(\hat x)-\alpha\varepsilon,
	\end{aligned}
	\end{equation*}
	thus, for $\overline{x}=\hat{x}+\eta^{'}$, we get
	\begin{equation}\label{equation_16}
	u(\overline{x})-v^*(\overline{x})\geq M-\alpha\varepsilon,
	\end{equation}
	in the case where $v^*(\hat{x})<\mathcal{H}_{inf}^\chi v^*(\hat{x})$ it suffice to choose $\eta^{'}\in V$ for which (\ref{equation_15}) holds, then we proceed analogously to get (\ref{equation_16}).\\
	Therefore, by taking $\alpha$ sufficiently small, we get $$\sup_{x\in\overline{B}_\delta(\overline{x})}\Bigl(u(x)-v^*(x)\Bigr)\geq M>0.$$
	As a consequence we consider $$M=\sup_{x\in\overline{B}_\delta(\overline{x})}\Bigl(u(x)-v^*(x)\Bigr).$$
	\par\textbf{Step 2.}
	Let $\varepsilon$ be a positive real number, $(x,y)\in\overline{B}_\delta(\overline{x})\times\overline{B}_\delta(\overline{x})$ and consider $\psi_\varepsilon$ the test function defined as follows:
	\begin{equation*}
	\begin{aligned}
		\psi_\varepsilon:\overline{B}_\delta(\overline{x})\times\overline{B}_\delta(\overline{x})&\rightarrow\mathbb{R}\\
		(x,y)&\rightarrow\psi_\varepsilon(x,y)=u(x)-v^*(y)-\frac{\parallel x-y\parallel ^2}{\varepsilon^2}.
	\end{aligned}
	\end{equation*}
	Let $(x_m,y_m)$ be the maximal point of $\psi_\varepsilon$ and denote $M_{\psi_\varepsilon}=\max\psi_\varepsilon(x,y)=\psi_\varepsilon(x_m,y_m).$ $M_{\psi_\varepsilon}$ exists, since on the one hand, $\psi_\varepsilon$ is a bounded and continuous function on a bounded set, and on the other hand, it is negative in a neighborhood of the boundary $\parallel x\parallel=\delta$ or $\parallel y\parallel=\delta$, and, by hypothesis, positive for some $x=y$. So, the search for the maximum can then be restricted to a compact set $B_{\delta-\gamma}(\overline{x})\times B_{\delta-\gamma}(\overline{x})$.\\
	Therefore
	\begin{equation}\label{equation_17}
	u(x_m)-v^*(y_m)-\frac{\parallel x_m-y_m\parallel ^2}{\varepsilon^2}\geq u(x)-v^*(y)-\frac{\parallel x-y\parallel ^2}{\varepsilon^2}.
	\end{equation}
	\begin{itemize}
		\item Firstly, for $x=x_m$, we get for all $y\in\overline{B}_\delta(\overline{x})$,
		$$u(x_m)-v^*(y_m)-\frac{\parallel x_m-y_m\parallel ^2}{\varepsilon^2}\geq u(x_m)-v^*(y)-\frac{\parallel x_m-y\parallel ^2}{\varepsilon^2},$$
		then $y_m$ is a local minimal point of $y\rightarrow (v^*-\phi_{v^*})(y)$ with
		$$\phi_{v^*}(y)=u(x_m)-\frac{\parallel x_m-y\parallel ^2}{\varepsilon^2}.$$
		From the inequality (\ref{equation_14}) and since $y_m\in\overline{B}_\delta(\overline{x})$ we get
		$$v(y_m)<\mathcal{H}_{inf}^\chi v(y_m),$$
		in addition $v$ is a viscosity supersolution to the new HJBI QVI, then by applying the result in Lemma \ref{lemma_4}, we find the following:
		$$min\Bigl[\lambda v^*(y_m)-D_y\phi_{v^*}(y_m).b(y_m)-f(y_m),\mathcal F_{inf}^c D_y\phi_{v^*}(y_m)\Bigr]> 0,$$
		thus, we obtain
		\begin{equation}\label{equation_18}
		\lambda v^*(y_m)-D_y\phi_{v^*}(y_m).b(y_m)-f(y_m)>0\;\text{and}\;\mathcal F_{inf}^c D_y\phi_{v^*}(y_m)>0.
		\end{equation}
		\item Secondly, for $y=y_m$, we get for all $x\in\overline{B}_\delta(\overline{x})$,
		$$u(x_m)-v^*(y_m)-\frac{\parallel x_m-y_m\parallel ^2}{\varepsilon^2}\geq u(x)-v^*(y_m)-\frac{\parallel x-y_m\parallel ^2}{\varepsilon^2},$$
		then $x_m$ is a local maximal point of $x\rightarrow (u-\phi_{u})(x)$ with
		$$\phi_{u}(x)=v^*(y_m)+\frac{\parallel x-y_m\parallel ^2}{\varepsilon^2}.$$
		Since $u$ is a viscosity subsolution to the new HJBI QVI, we get $$u(x_m)-\mathcal{H}_{inf}^\chi u(x_m)\leq 0,$$
		and
		\begin{equation}\label{equation_19}
		\lambda u(x_m)-D_x\phi_{u}(x_m).b(x_m)-f(x_m)\leq 0\;\text{or}\;\mathcal F_{inf}^c D_x\phi_{u}(x_m)\leq 0.
		\end{equation}
	\end{itemize}
	From (\ref{equation_18}) and (\ref{equation_19}), since $\mathcal F_{inf}^c D_y\phi_{v^*}(y_m)=\mathcal F_{inf}^c D_x\phi_{u}(x_m)$, we get
	$$\lambda v^*(y_m)-D_y\phi_{v^*}(y_m).b(y_m)-f(y_m)>0$$
	and
	$$\lambda u(x_m)-D_x\phi_{u}(x_m).b(x_m)-f(x_m)\leq 0.$$
	It follows that
	$$\lambda\Bigl(u(x_m)-v^*(y_m)\Bigr) +\Bigl[D_y\phi_{v^*}(y_m).b(y_m)-D_x\phi_{u}(x_m).b(x_m)\Bigr]+\Bigl(f(y_m)-f(x_m)\Bigr)\leq 0.$$
	Then
	$$\lambda\Bigl(u(x_m)-v^*(y_m)\Bigr) +2\frac{\parallel x_m-y_m\parallel }{\varepsilon^2}\Bigl(b(y_m)-b(x_m)\Bigr)+\Bigl(f(y_m)-f(x_m)\Bigr)\leq 0.$$
	Since $f$ and $b$ are Lipschitz with constants $C_f$ and $C_b$, respectively, we deduce
	\begin{equation}\label{equation_20}
	\lambda\Bigl(u(x_m)-v^*(y_m)\Bigr)\leq 2C_b\frac{\parallel x_m-y_m\parallel ^2}{\varepsilon^2}+C_f\parallel x_m-y_m\parallel.
	\end{equation}
	\par\textbf{Step 3.}
	Proving now that $\forall\beta>0,\;\exists\varepsilon_0>0,\;\forall\varepsilon\leq\varepsilon_0:\;\frac{\parallel x_m-y_m\parallel ^2}{\varepsilon^2}\leq\beta$ and showing the contradiction.\\
	By taking $x=y$ in (\ref{equation_17}) we get, for all $x\in\overline{B}_\delta(\overline{x})$, $u(x)-v^*(x)\leq M_{\psi_\varepsilon}$ then $0<M\leq M_{\psi_\varepsilon}$.\\
	Let $r^2=\parallel u\parallel_\infty+\parallel v^*\parallel_\infty$ then $0<M_{\psi_\varepsilon}\leq r^2-\frac{\parallel x_m-y_m\parallel ^2}{\varepsilon^2}$, it follows that $\parallel x_m-y_m\parallel\leq\varepsilon r$.\\
	Since $u$ is upper semi-continuous, we have
	$$\forall\beta>0,\;\exists\varepsilon_0>0,\;\forall\varepsilon\leq\varepsilon_0:\;u(x_m)\leq u(y_m)+\beta,$$
	then
	$$u(x_m)-v^*(y_m)\leq u(y_m)-v^*(y_m)+\beta\leq M+\beta.$$
	Using $M\leq M_{\psi_\varepsilon}$ we get
	$$\frac{\parallel x_m-y_m\parallel^2}{\varepsilon^2}\leq\beta\;\text{and}\;M\leq u(x_m)-v^*(y_m),$$
	then from (\ref{equation_20}) we deduce
	$$\lambda\Bigl(u(x_m)-v^*(y_m)\Bigr)\leq2C_b\beta+C_f\varepsilon\sqrt{\beta}.$$
	By sending $\beta$ to 0, we get $u(x_m)-v^*(y_m)\leq 0,$ from which yields for all $x\in\overline{B}_\delta(\overline{x})$,
	$$u(x)-v^*(x)\leq u(x_m)-v^*(y_m)\leq 0,$$
	thus we get the contradiction $M\leq 0$.\\
	Hence, by letting $\mu\rightarrow 1$ and $\alpha\rightarrow 0$, we deduce for all $x\in\mathbb{R}^n$ the desired comparison:
	$$u(x)\leq v(x).$$
\end{proof}

\begin{theorem}
	Under assumptions \textbf{H1} and \textbf{H2}, the new HJBI QVI has a unique bounded and continuous viscosity solution.
\end{theorem}

\begin{proof}
	Assume that $v_1$ and $v_2$ are two viscosity solutions to the new HJBI QVI. We first use $v_1$ as a bounded and continuous viscosity subsolution and $v_2$ as a bounded and continuous viscosity supersolution and we recall the comparison theorem. Then we change the role of $v_1$ and $v_2$ to get $v_1(x)=v_2(x)$ for all $x\in\mathbb{R}^n$.
\end{proof}

\begin{corollary}
	Under assumptions \textbf{H1} and \textbf{H2}, the lower and upper value functions coincide and the value function $V:=V^-=V^+$ of the infinite horizon two-player zero-sum deterministic differential game is the unique viscosity solution to the new HJBI QVI.
\end{corollary}

\end{document}